\documentclass[11pt]{article}
\usepackage{amsmath}
\usepackage{amsfonts}
\usepackage{graphicx}
\usepackage{latexsym}
\usepackage{color}

\newcommand{\R}{\mathbb{R}}
\newcommand{\E}{\mathbb{E}}
\newcommand{\PP}{\mathbb{P}}

\newcommand{\N}{\mathbb{N}}

\newcommand{\De}{\Delta}

\newcommand{\stab}{\stackrel{d_{st}}{\longrightarrow}}
\newcommand{\schw}{\stackrel{d}{\longrightarrow}}
\newcommand{\toop}{\stackrel{\PP}{\longrightarrow}}

\newcommand{\ucp}{\stackrel{\mbox{\tiny u.c.p.}}{\Longrightarrow}}

\newcommand{\bee}{\begin{equation}}
\newcommand{\eee}{\end{equation}}
\newcommand{\bea}{\begin{eqnarray}}
\newcommand{\eea}{\end{eqnarray}}
\newcommand{\bean}{\begin{eqnarray*}}
\newcommand{\eean}{\end{eqnarray*}}

\newcommand{\qed}{$\hfill\Box$}

\newtheorem{prop}{Proposition}[section]

\newtheorem{lem}[prop]{Lemma}
\newtheorem{ex}[prop]{Example}
\newtheorem{theo}[prop]{Theorem}
\newtheorem{rem}[prop]{Remark}

\begin{document}

\title{On non-standard limits of Brownian semi-stationary processes}
\author{Kerstin G\"artner \thanks{Institute of Population Genetics, University of Veterinary Medicine,
Veterin\"arplatz 1, 1210 Vienna, Austria,
Email: kerstin.gaertner@vetmeduni.ac.at.} \and
Mark Podolskij \thanks{Department of Mathematics, University of Aarhus,
Ny Munkegade 118, 8000 Aarhus C,
Denmark, Email: mpodolskij@creates.au.dk.}}

\date{\today}

\maketitle

\begin{abstract}
In this paper we present some new asymptotic results for high frequency statistics
of Brownian semi-stationary ($\mathcal{BSS}$) processes. More precisely, we will show
that singularities in the weight function, which is one of the ingredients of a $\mathcal{BSS}$
process, may lead to non-standard limits of the realised quadratic variation. In this case
the limiting process is a convex combination of shifted integrals of the intermittency function.
Furthermore, we will demonstrate the corresponding stable central limit theorem. Finally, 
we apply the probabilistic theory to study the asymptotic properties of the realised ratio statistics,
which estimates the smoothness parameter of a $\mathcal{BSS}$
process.

\ \

\bigskip


\end{abstract}

\section{Introduction}
\label{Intro}
\setcounter{equation}{0}
\renewcommand{\theequation}{\thesection.\arabic{equation}}
In the last years Brownian semi-stationary processes and their tempo-spatial extensions, \textit{ambit fields}, have
been widely studied in the literature. This class of models has been originally proposed by Barndorff-Nielsen and Schmiegel
\cite{BS07} in the context of turbulence modeling. In their general form, Brownian semi-stationary processes without drift 
are defined as
\begin{equation*}
X_{t}=\mu +\int_{-\infty }^{t}g(t-s)\sigma _{s}W(\mathrm{d}s), \qquad t\in \R
\end{equation*}%
where $\mu $ is a constant, $W$ is a Brownian measure on $\mathbb{R}$, $g:\R\rightarrow \R$
is a deterministic weight function with
$g( t) =0$ for $t\leq 0$, and $\sigma $ is a c\`{a}dl\`{a}g processes. If  $\sigma$ is stationary and independent of $W$, 
then $(X_t)_{t\in \R}$ is stationary, which explains the name Brownian semi-stationary process. In the framework
of turbulence modeling, $(X_t)_{t\in \R}$ denotes the velocity of a turbulent flow in the direction of the mean field
measured at a fixed point in space. The stochastic process $(\sigma_t)_{t\in \R} $
embodies the \emph{intermittency} of the dynamics of $X$. We refer to \cite{BS07, BS08a, BS08b, BS09} for application
of Brownian semi-stationary processes and ambit fields to turbulence modeling, and to \cite{BBC, BJJS07} for further
applications in mathematical finance and biology.

Recently, probabilistic properties of high frequency statistics of $\mathcal{BSS}$ processes have been investigated 
in several papers. We refer to a series of articles \cite{BCP11,BCP13,CHPP13}, which studies the asymptotic behaviour
of (multi)power variation of $\mathcal{BSS}$ models. Typically, the weight function $g$ considered in the aforementioned 
work has the form
\[
g(x)=x^{\alpha} f(x),\qquad \alpha \in (-1/2,0) \cup (0,1/2),
\]  
where $f$ is a sufficiently smooth function slowly varying at $0$ and with rapid decay at infinity. This type 
of weight functions satisfies $g\in \mathbb L^2 (\R)$, but $g' \not \in \mathbb L^2 (\R)$ since $g'$ is not square integrable
near $0$; in other words, the latter property means that $0$ is the only \textit{singularity point} of the 
weight function $g$. As a consequence, the process $X$ is not a semimartingale. Moreover, its local behaviour corresponds to the 
one of a fractional Brownian motion with Hurst parameter $H=\alpha +1/2$. 

Understanding the limit theory for $\mathcal{BSS}$ processes requires an analysis of the following probability measure.
For any $A\in \mathcal{B} (\R)$, we define
\begin{align} \label{pindef}
\pi_n (A) := \frac{\int_A \{g(x+\Delta_n) -g(x)\}^2 dx}{\int_{\R} \{g(x+\Delta_n) -g(x)\}^2 dx}.
\end{align}
In the setting of weight functions as above, we deduce that $\pi_n \schw \delta_0$ as $\Delta_n\rightarrow 0$, 
where $\delta_0$ denotes the Dirac measure
at $0$ (cf. \cite{BCP11}). In this case the limit of the power variation of a $\mathcal{BSS}$ process is given as
\begin{align} \label{standarducp} 
\Delta_n \tau_n^{-p}\sum_{i=1}^{[t/\Delta_n]} |X_{i\Delta_n} - X_{(i-1)\Delta_n}|^p  \ucp m_p \int_0^t |\sigma_s|^p ds, 
\qquad \text{as } \Delta_n\rightarrow 0,
 \end{align} 
where $m_p=\E[|\mathcal{N}(0,1)|^p]$, $\tau_n$ is a certain normalizing sequence and $\ucp$ stands
for convergence in probability uniformly on compact sets. In \cite{BCP11,BCP13} the asymptotic mixed
normality of (multi)power variation is proved and the paper \cite{CHPP13} studies the application of the limit theory
to estimation of the \textit{smoothness parameter} $\alpha$. We remark that the asymptotic results are quite similar to the
theory of power variations of continuous It\^o semimartingales (cf. \cite{BGJPS06,J08} among many others), although the methodologies
of proofs are completely different.  

The aim of this paper is to demonstrate that other type of limits for power variations may appear when 
the weight function $g$ exhibits further singularity points. More precisely, we will prove that 
\begin{align} \label{nonstandarducp}
\Delta_n \tau_n^{-2}\sum_{i=1}^{[t/\Delta_n]} (X_{i\Delta_n} - X_{(i-1)\Delta_n})^2  \ucp  
\int_0^\infty \left( \int_{-\theta}^{t-\theta } \sigma_s^2 ds\right) \pi (d\theta), 
\qquad \text{as } \Delta_n\rightarrow 0,
 \end{align} 
where $\pi$ is a probability measure on $[0,\infty )$ whose support is a subset of all singularity points of $g$. 
Consequently, the limit theory for $\mathcal{BSS}$ processes is richer than the corresponding theory for continuous
It\^o semimartingales.  
Furthermore,
we will show the associated stable central limit theorem including the setting of higher order differences. We remark that 
this type of limits has already appeared in \cite{BS07}. The authors proved convergence in probability for the realised 
quadratic variation under the independence assumption between $\sigma$ and $W$, and under further conditions on certain measures
associated with $g$, which identify $\pi$. However, it remained quite unclear when a given weight function $g$ satisfies 
the proposed set of conditions. The main goal of our paper is to show that singularity points of $g$, i.e. all points 
around which $g'$ is not square integrable, determine the support and the weights of $\pi$. 
We remark that in physics multiple singularity points of $g$ lead to non-homogeneous turbulent flows. 
Moreover, we will study 
the effect of this new class of weight functions $g$ on smoothness parameter estimation. More precisely, we will present
the asymptotic behaviour of a realised ratio statistic that compares the realised quadratic variation at two different 
frequencies.

The paper is organised as follows. Section 2 presents the main framework and a set of assumptions. In Section 3
we demonstrate the complete asymptotic theory for the realised quadratic variation of $\mathcal{BSS}$ processes,
including the law of large numbers and the associated stable central limit theorem. In Section 4 we apply the probabilistic
results to determine the asymptotic behaviour of a realised ratio statistic, which is an estimator of the smoothness
parameter of $X$. Finally, all proofs are collected in Section 5.

\section{The setting} \label{sec2}
\setcounter{equation}{0}
\renewcommand{\theequation}{\thesection.\arabic{equation}}

\subsection{Model} \label{sec2.1}
We start with a given filtered probability space  $(\Omega ,\mathcal{F},(\mathcal{F}_{t})_{t\in \mathbb{R}},\mathbb{P})$
on which our processes are defined. We consider a $\mathcal{BSS}$ process $(X_{t})_{t\in \mathbb{%
R}}$ (without drift) given as
\begin{equation} \label{x}
X_{t}=\mu +\int_{-\infty }^{t}g(t-s)\sigma _{s}W(ds),\qquad t\in \mathbb{R},
\end{equation}
where $W$ is an $(\mathcal{F}_{t})_{t\in \mathbb{R}}$-adapted white noise on $\mathbb{R}$, 
$g:\R \rightarrow \mathbb{R}$ is a deterministic weight function satisfying
$g(t)=0$ for $t\leq 0$ and $g\in \mathbb{L}^{2}(\mathbb{R})$. The intermittency  process $\sigma $ 
is assumed to be an \ $(\mathcal{F}_{t})_{t\in \mathbb{R}}$-adapted c\`{a}dl\`{a}g process. We recall that 
$\left( \mathcal{F}_{t}\right) $-adapted white noise $W$ is zero-mean Gaussian random measure on $\{A\in \mathcal{B}(\R):~
\lambda (A)<\infty \}$, where $\lambda$ denotes the Lebesgue measure, whose covariance kernel is given by
\[
\E[W(A)W(B)] = \lambda (A\cap B).
\] 
The finiteness of the process $X$ is guaranteed by the condition
\begin{align} \label{fincon}
\int_{-\infty }^{t}g^2(t-s)\sigma^2_{s} ds<\infty \quad \text{almost surely}, 
\end{align}
for any $t\in \R$, which we assume from now on. The underlying observations of the $\mathcal{BSS}$ process $X$ are
\[
X_0, X_{\Delta_n}, X_{2\Delta_n}, \ldots, X_{\Delta_n [t/\Delta_n]}
\] 
with $\Delta_n\rightarrow 0$ and $t$ fixed. In other words, we are in the framework of infill asymptotics. Our realised
quadratic variation statistics will be based upon higher order increments of $X$ computed at different frequencies.
For any $k\in \N$ and $v=1,2$, the $k$-th order difference $\Delta_{i,k}^{n,v} X$ at frequency $v\Delta_n$ and at stage $i\geq vk$ is defined by
\begin{align} \label{filter}
\Delta_{i,k}^{n,v} X:= \sum_{j=0}^k (-1)^j \binom{k}{j} X_{(i-vj)\Delta_n}.
\end{align}
The quantity $\Delta_{i,k}^{n,v} X$ is a particular example of a $k$-th order filter applied to the process $X$. 
When $v = 1$ we usually write $\Delta_{i,k}^{n}X$ instead of $\Delta_{i,k}^{n,1}X$. For instance, 
\[
\Delta_{i,1}^n X = X_{i\Delta_n}-X_{(i-1)\Delta_n} \quad \textrm{and} \quad \Delta_{i,2}^n X = X_{i\Delta_n}-2X_{(i-1)\Delta_n}+X_{(i-2)\Delta_n}.
\] 
The realised quadratic variation statistic based upon $\Delta_{i,k}^{n,v} X$ is defined as
\begin{align} \label{qv}
QV(X,k,v\Delta_n)_t &:= \sum_{i=vk}^{[t/\Delta_n]} (\Delta_{i,k}^{n,v} X)^2
\end{align}
As in \cite{BCP11,BCP13}, the \textit{Gaussian core} $G$ is crucial for understanding the fine structure of $X$. The process
$G=(G_t)_{t\in \R}$ is a zero-mean stationary Gaussian process given by
\begin{align} \label{g}
G_t: = \int_{-\infty}^t g(t-s)  W(ds), \qquad t\in \R.
\end{align}  
We remark that $G_t<\infty$ since $g\in \mathbb{L}^2(\R)$. A straightforward computation
shows that the correlation kernel $r$ of $G$ has the form
\[
r(t) = \frac{\int_0^\infty g(u)g(u+t) du}{\|g \|_{\mathbb L^2(\R)}^2}, \qquad t\geq 0.
\]
Another important quantity for the asymptotic theory is the variogram  $R$, i.e.
\begin{align} \label{R}
R(t):= \E[(G_{t+s} - G_s)^2] = 2 \|g \|_{\mathbb L^2(\R)}^2 (1-r(t)), \qquad 
\tau_{k}(v\Delta_n):=\sqrt{\E[(\Delta_{i,k}^{n,v} G)^2]}.
\end{align} 
The quantity $\tau_{k}(v\Delta_n)$ will appear as a proper scaling in the law of large numbers
for the statistic $QV(X,k,v\Delta_n)$.

\subsection{Main assumptions} \label{sec2.2}
As we mentioned in the introduction, understanding the asymptotic behaviour of the probability measure 
\[
\pi_n (A) = \frac{\int_A \{g(x+\Delta_n) -g(x)\}^2 dx}{\int_{\R} \{g(x+\Delta_n) -g(x)\}^2 dx}, \qquad A\in \mathcal{B}(\R),
\]
is absolutely crucial for determining the limit theory for the realised quadratic variation $QV(X,1,v\Delta_n)$
(for $QV(X,k,v\Delta_n)$, $k\geq 2$, there exists an analogous probability measure). Indeed, the condition 
\[
\pi_n \schw \pi, 
\]  
where $\pi$ is a certain probability measure on $\R$, is necessary (but not sufficient) to obtain a non-standard law 
of large numbers at \eqref{nonstandarducp}. In \cite{BCP11,BCP13,CHPP13} it has been dealt with the case $\pi =\delta_{0}$,
and hence we obtained a rather standard convergence as in \eqref{standarducp}. However, due to a moving
average structure of the process $X$, even trivial weight functions $g$
may lead to  $\pi \not=\delta_{0}$ as the next simple example shows.

\begin{ex} \label{ex1} \rm
Let us consider the function $g(x)=1_{[0,1]}(x)$. A simple computation shows that 
\[
\pi_n (A)= \frac{\{\lambda (A\cap [-\Delta_n,0]) + \lambda (A\cap [1-\Delta_n,1]) \}}{2\Delta_n},
\]
and consequently $\pi_n \schw \pi=\frac{1}{2}(\delta_0+\delta_1)$. Indeed, the convergence in \eqref{nonstandarducp}
with $\tau_n^2=\tau_{1}(\Delta_n)^2=2\Delta_n$ can be shown in a straightforward manner. For our weight function $g$, we deduce that
\[
X_t=Y_t-Y_{t-1} \qquad \text{with} \qquad Y_t=\int_{-1}^t \sigma_s W(ds),
\] 
for $t\geq -1$. Noticing that $Y$ is a martingale, we easily conclude   
\[
\Delta_n \tau_n^{-2}\sum_{i=1}^{[t/\Delta_n]} (X_{i\Delta_n} - X_{(i-1)\Delta_n})^2  \ucp \frac 12
\left( \int_{0}^{t} \sigma_s^2 ds + \int_{-1}^{t-1} \sigma_s^2 ds\right), 
\qquad \text{as } \Delta_n\rightarrow 0,
\]
which confirms \eqref{nonstandarducp}. This example demonstrates that if $g(x)=\sum_{i=1}^{l}a_i1_{[\theta_{i}^{(1)},
\theta_i^{(2)}]}$ 
with $0\leq \theta_1^{(1)}<\theta_1^{(2)}<\theta_2^{(1)}<\theta_2^{(2)}<\cdots<\theta_l^{(2)}<\infty$ then
\[
\text{supp}(\pi)=\{\theta_1^{(1)}, \theta_1^{(2)}, \ldots, \theta_l^{(1)},\theta_l^{(2)}\} \qquad \text{and} \qquad 
\pi (\{ \theta_i^{(1)}\}) = \pi (\{ \theta_i^{(2)}\})= \frac{a_i^2}{2\sum_{i=1}^{l}a_i^2}, 
\] 
and \eqref{nonstandarducp} holds.
\qed
\end{ex}
Barndorff-Nielsen and Schmiegel \cite{BS07} provide conditions on certain rather complex measures associated with 
$g$ (including $\pi_n \schw \pi$), which are sufficient for proving law of large numbers of the type \eqref{nonstandarducp}
under the independence assumption between $\sigma$ and $W$. However, it is not a priori clear when a given weight function 
$g$ satisfies those conditions. Furthermore, conditions ensuring the associated central limit theorem are expected to be even
more complex. 

In this paper we follow a different route. We present  an explicit large class of weight functions $g$, which leads to 
the law of large numbers of \eqref{nonstandarducp}, such that the limiting probability measure $\pi$ is easily identified.
Moreover, the associated central limit theorems are obtained (the limit theory does not require independence
of $\sigma$ and $W$). The crucial message of this paper is that singularity points of $g$ defined below determine
the support and the weights of $\pi$.   

Let $0=\theta_0<\theta_1<\cdots<\theta_l<\infty $ be a set of given points and $\alpha_0,\ldots ,\alpha_l\in (-1/2,0)\cup (0,1/2)$.
For any function $h\in C^m (\R)$, $h^{(m)}$ denotes the $m$-th derivative of $h$. Recall that $k\geq 1$ stands for the order of the
filter defined in \eqref{filter}. We introduce the following set of assumptions. \\ \\
(A): For $\delta<\frac 12 \min_{1\leq i\leq l}(\theta_i-\theta_{i-1})$ it holds that \\ \\
(i) $g(x)= x^{\alpha_0}  f_0(x)$ for $x\in (0, \delta)$ and $g(x)= |x-\theta_l|^{\alpha_l}  f_l(x)$
for $x\in (\theta_l-\delta, \theta_l)\cup (\theta_l, \infty )$. \\ \\
(ii) $g(x)= |x-\theta_i|^{\alpha_i}  f_i(x)$ for $x\in (\theta_i-\delta, \theta_i)\cup (\theta_i, \theta_i +\delta)$,
$i=1,\ldots,l-1$. \\ \\
(iii) $g(\theta_i)=0$, 
$f_i\in C^k\left( (\theta_i-\delta, \theta_i +\delta)\right)$ and $f_i(\theta_i)\not =0$ for $i=0,\ldots, l$.\\ \\
(iv) $g\in C^k(\R \setminus \{\theta_0, \ldots, \theta_l\})$ and $g^{(k)}\in \mathbb L^2 \left(
\R \setminus \cup_{i=0}^l (\theta_i-\delta, \theta_i+\delta) \right)$. \\ \\
(v) For any $t>0$
\begin{align} \label{F}
F_t= \int_{\theta_l+1}^\infty g^{(k)}(s)^2 \sigma_{t-s}^2 ds <\infty.  
\end{align} 
We also set
\begin{align} \label{min}
\alpha := \min\{\alpha_0,\ldots ,\alpha_l \}, \qquad \mathcal{A}:=\{0\leq i\leq l:~\alpha_i=\alpha\}.
\end{align}   
Let us give some remarks on this set of conditions.

\begin{rem} \label{rem1} \rm
The points $\theta_0,\ldots,\theta_l$ are singularities of $g$ in the sense that $g^{(k)}$ is not square integrable 
around these points, because $\alpha_0,\ldots ,\alpha_l\in (-1/2,0)\cup (0,1/2)$ and conditions (A)(i)-(iii) hold. 
Condition (A)(iv) indicates that $g$ exhibits no further singularities. The papers \cite{BCP11,BCP13,CHPP13} 
deal with the framework of a single singularity at $0$. \qed
\end{rem}

\begin{rem} \label{rem2} \rm
The parameter $\alpha \in (-1/2,0)\cup (0,1/2)$ defined at \eqref{min} determines the smoothness coefficient
of the $\mathcal{BSS}$ process $X$. In some sense, the coefficients $\alpha_i$ with $i\in \mathcal{A}$
will dominate when proving the limit theory for $QV(X,k,v\Delta_n)_t$. In particular, we will prove
that $\text{supp}(\pi)=\{\theta_i\}_{i\in \mathcal{A}}$.   \qed
\end{rem}

\begin{rem} \label{rem3} \rm
The weight function considered in Example \ref{ex1} obviously does not satisfy the assumption (A).  
Indeed, in the framework of Example \ref{ex1}
the limit theory for $QV(X,k,v\Delta_n)_t$ relies on semimartingale methods (cf. \cite{BGJPS06})
as $X$ is a difference of two martingales (although $X$ is not a semimartingale). In the case
of assumption (A) we are in the framework of fractional processes. More precisely, the small scale
behaviour of the Gaussian core $G$ of $X$ is close to the small scale bahaviour of a fractional Brownian
motion with Hurst parameter $H=\alpha + 1/2$. In this situation the limit theory for $QV(X,k,v\Delta_n)_t$
relies on Malliavin calculus and Bernstein's blocking technique.   \qed
\end{rem}

\begin{rem} \rm
In papers \cite{BCP11,BCP13,CHPP13}, where $l=0$ holds, the function $f_0$ is assumed to be slowly varying at $\theta_0=0$.
In this setting more assumptions are required to establish the limit theory than mere condition (A). In our paper
we impose a bit stronger assumptions on functions $f_j$, $j=0,\ldots,l$, in order to avoid a longer set of further
conditions.

Note that condition (A)(ii) implies a symmetric behaviour of the function $g$ around the points $\theta_j$, $j=1,\ldots,l$.
Instead we could have assumed different power behaviour left and right from $\theta_j$. Although certain constants
in the limit theorems would change in this case, the asymptotic theory remains essentially the same. \qed  
\end{rem}

\section{Limit theorems} \label{sec3}
\setcounter{equation}{0}
\renewcommand{\theequation}{\thesection.\arabic{equation}}

\subsection{Law of large numbers} \label{sec3.1}
For any number $k\geq 1$ and $v=1,2$, we introduce a $k$-th order filter associated with $g$ via
\begin{align} \label{filterg}
\Delta_{k}^{n,v} g(x):= \sum_{j=0}^k (-1)^j \binom{k}{j} g(x-vj\Delta_n), \qquad x\in \R.
\end{align}
There is a straightforward relationship between the scaling quantity $\tau_{k}(v\Delta_n)$
defined at \eqref{R} and the function $\Delta_{k}^{n,v} g$, namely
\[
\tau_{k}(v\Delta_n)^2 = \|\Delta_{k}^{n,v} g\|_{\mathbb{L}^2(\R)}^2.
\] 
Now, we define the corresponding measures associated with $\Delta_{k}^{n,v} g$:
\begin{align} \label{pink}
\pi_{n,k}^{v}(A) := \frac{\int_A (\Delta_{k}^{n,v} g(x))^2 dx}{\|\Delta_{k}^{n,v} g\|_{\mathbb{L}^2(\R)}^2},
\qquad A\in \mathcal{B}(\R).
\end{align}
In order to identify the limit of $\pi_{n,k}^{v}$, we define the following functions
\begin{align} \label{hdef}
h_0 (x)&:= f_0(\theta_0) \sum_{j=0}^k (-1)^j \binom{k}{j} (x-j)_{+}^{\alpha_0}, \\[1.5 ex]
h_i (x)&:= f_i(\theta_i) \sum_{j=0}^k (-1)^j \binom{k}{j} |x-j|^{\alpha_i}, \qquad i=1,\ldots, l, \nonumber
\end{align}
where $x_+:=\max\{x,0\}$. At this stage we suppress the dependence of functions $h_i$ on the index $k$.
Our first result presents the limiting measure $\pi_k$, which will appear in the law of large numbers.

\begin{prop} \label{prop1}
Assume that condition (A) holds. Then we deduce that
\begin{align*}
\pi_{n,k}^{v} \schw \pi_k,
\end{align*}
for any $k\geq 1$ and $v=1,2$, where the probability measure $\pi_k$ is given as
\begin{align} \label{pik}
\text{supp}(\pi_k)= \{\theta_i\}_{i\in \mathcal{A}}, \qquad 
\pi_k (\theta_i)= \frac{\|h_i\|_{\mathbb{L}^2(\R)}^2 1_{i\in \mathcal{A}}}{\sum_{j=0}^l \|h_j\|_{\mathbb{L}^2(\R)}^2
1_{j\in \mathcal{A}}}.
\end{align}
\end{prop}
Recalling the definition of the set $\mathcal{A}$ at \eqref{min}, Proposition \ref{prop1} says that only
singularities corresponding to the minimal indexes $\alpha_i$ (i.e. indexes with $\alpha_i=\alpha$)
contribute to the limit. We remark that the norms $\|h_i\|_{\mathbb{L}^2(\R)}$ are indeed finite, since
for $|x|$ large enough 
\[
|h_i(x)|^2\leq C |x|^{2(\alpha_i-k)} \qquad \text{and} \qquad  2(\alpha_i-k)<-1,
\]
for any $k\geq 1$ and $\alpha_i\in (-1/2,0)\cup (0,1/2)$ due to Taylor expansion.
Our next result is the law of large numbers for the statistic
$QV(X,k,v\Delta_n)$.

\begin{theo} \label{th1}
Assume that condition (A) holds. Then 
\begin{align} \label{LLN}
\frac{\Delta_n}{\tau_{k}(v\Delta_n)^2} QV(X,k,v\Delta_n)_t \ucp   QV(X,k)_t:= 
\int_0^\infty \left( \int_{-\theta}^{t-\theta } \sigma_s^2 ds\right) \pi_k (d\theta),
\end{align}
where the probability measure $\pi_k$ is introduced in \eqref{pik}.
\end{theo} 


\subsection{Central limit theorem}
Now, we will present a stable central limit theorem associated with convergence in \eqref{LLN}. Let us shortly 
recall the notion of stable convergence, which is originally due to R\'enyi \cite{R63}.  
We say that a sequence of processes $Y^n$ converges stably in law to
a  process $Y$, where $Y$ is defined on an extension $(\Omega', \mathcal{F}', \mathbb P')$ 
of the original probability space $(\Omega, \mathcal{F}, \mathbb P)$, in the space $\mathbb D([0,T])$
equipped with the uniform topology ($Y^n \stab Y$) if and only if
\begin{equation*} 
\lim_{n\rightarrow \infty} \E[f(Y^n) Z] = \E'[f(Y)Z] 
\end{equation*} 
for any bounded and continuous function $f: \mathbb D([0,T]) \rightarrow \R$ and any bounded $\mathcal F$-measurable
random variable $Z$.  We refer to \cite{AE78,R63} for a detailed study of stable convergence. 
Note that stable convergence is a stronger mode of convergence than weak convergence, but it is weaker than 
u.c.p. convergence.

The stable central limit theorem associated with convergence in \eqref{LLN} is different compared to the corresponding
result in the case of a single singularity (cf. \cite{BCP11,BCP13}). In particular, as we will see below, the limiting process
is not an $\mathcal F$-conditional martingale on every interval $[0,T]$, but just for small enough $T$. 
For the purpose of statistical inference we present a joint central limit theorem for the pair 
$(QV(X,k,\Delta_n),QV(X,k,2\Delta_n))$.

\begin{theo} \label{th2}
Assume that condition (A) holds and the intermittency process $\sigma$ is H\"older continuous of order
$\gamma >1/2$. If $k=1$ we further assume that $\alpha_j \in (-\frac 12,0)$ for all $0\leq j\leq l$. Then, under condition
\begin{align} \label{robustness}
\alpha_i -\alpha > 1/4 \qquad \text{for all } i \not \in \mathcal{A},  
\end{align} 
we obtain the stable convergence
\begin{align} 
&\Delta_n^{-1/2} \left(\frac{\Delta_n}{\tau_{k}(\Delta_n)^2} QV(X,k,\Delta_n)_t -   QV(X,k)_t,
\frac{\Delta_n}{\tau_{k}(2\Delta_n)^2} QV(X,k,2\Delta_n)_t -   QV(X,k)_t \right)^{\star} \nonumber \\[1.5 ex] 
& \label{CLT} \stab  L_t= \int_0^t \mathcal V_s^{1/2} dB_s
\end{align}
on $\mathbb D^2([0, \min_{1\leq j\leq l}(\theta_j-\theta_{j-1})])$ equipped with the uniform topology, where 
$B$ is a two-dimensional Brownian motion, independent of $\mathcal F$, defined on an
extension of the original probability space $(\Omega, \mathcal{F}, \mathbb P)$. The matrix $\mathcal V_s$ is given by 
\begin{align} \label{asyvar}
\mathcal V_s= \left(\int_0^{\infty} \sigma^2_{s-\theta} \pi_k (d\theta)\right)^2 \Lambda_k,
\end{align}
where the $2\times 2$ matrix $\Lambda_k=(\lambda_{ij}^k)_{1\leq i,j\leq 2}$ is  defined by
\begin{align} 
\lambda_{11}^k &= \lim_{n\rightarrow \infty } \De_n^{-1} \mathrm{var} \Big( 
\frac{\Delta_n}{\hat{\tau}_{k}(\Delta_n)^2} QV(B^H,k,\Delta_n)_1\Big), \nonumber\\
\label{lambda} \lambda_{22}^k &= \lim_{n\rightarrow \infty } \De_n^{-1} \mathrm{var} \Big( 
\frac{\Delta_n}{\hat{\tau}_{k}(2\Delta_n)^2} QV(B^H,k,2\Delta_n)_1\Big)  \\
\lambda_{12}^k &= \lim_{n\rightarrow \infty } \De_n^{-1} \mathrm{cov} \Big( 
\frac{\Delta_n}{\hat{\tau}_{k}(\Delta_n)^2} QV(B^H,k,\Delta_n)_1, 
\frac{\Delta_n}{\hat{\tau}_{k}(2\Delta_n)^2} QV(B^H,k,2\Delta_n)_1 \Big) \nonumber
\end{align}
with $B^H$ being a fractional Brownian motion with Hurst parameter $H=\alpha + 1/2$ and 
$\hat{\tau}_{k}(v\Delta_n)^2:=\E[(\Delta_{i,k}^{n,v} B^H)^2]$. 
\end{theo}
The H\"older condition  is a standard requirement for the validity of the blocking technique 
applied in the proofs (cf. \cite{BCP11,BCP13}). As we remarked earlier, the singularity points $\theta_i$ with
$i \not \in \mathcal{A}$ do not affect the law of large numbers in \eqref{LLN}. However, they are responsible
for a certain bias, which might explode in the central limit theorem. Assumption \eqref{robustness} guarantees
that it does not happen.  

The appearance of the fractional Brownian motion in the definition of the matrix $\Lambda_k$ is explained
by the fact that the local behaviour of the Gaussian core $G$ is close to the local behaviour of $B^H$
with $H=\alpha +1/2$. In the terminology of the theory of Gaussian fields it means that $B^H$
is a \textit{tangent process} of $G$. In particular, the correlation structure of increments of $G$
converges to the  correlation structure of increments of $B^H$. 

\begin{rem} \label{rem6} \rm
The limiting process $L$ is an $\mathcal F$-conditional Gaussian martingale on the interval 
$[0, \min_{1\leq j\leq l}(\theta_j-\theta_{j-1})]$. Outside of this interval the $\mathcal F$-conditional martingale
property gets lost. One may still show a stable central limit theorem with an $\mathcal F$-conditional Gaussian
process as the limit, but only when $\theta_j-\theta_{j-1}\in \mathbb N$ for all $j$, since otherwise the covariance structure
of the original statistic does not converge. We dispense with the exact presentation of this case.
\qed
\end{rem} 

\begin{rem} \label{rem4} \rm
The limits in \eqref{lambda} are indeed finite and can be computed explicitly. To see this, let us define the fractional Brownian noise
of order $k$ and scale $v=1,2$ via
\begin{align} \label{fbmnoise}
\Delta_{i,k}^v B^H:= \sum_{j=0}^k (-1)^j \binom{k}{i} B^H_{i-vj}, 
\end{align}
and set 
\begin{align} \label{fbmcorr}
\rho_k^{v_1,v_2}(j):=\text{corr}(\Delta_{i,k}^{v_1} B^H, \Delta_{i+j,k}^{v_2} B^H)
\end{align}
(Recall that $B^H$ has stationary increments.)
Using the covariance kernel of the fractional Brownian motion one can compute the quantity
$\rho_k^{v_1,v_2}(j)$ explicitly. For instance,
\[
\rho_1^{1,1}(j)= \frac{1}{2} \Big(|j+1|^{2H} - 2|j|^{2H} + |j-1|^{2H} \Big), \qquad j\geq 1.
\] 
A straightforward computation shows that $|\rho_k^{v_1,v_2} (j)|= O( |j|^{2(H-k)})$ as $|j|\rightarrow \infty$. Hence, using $H$-self
similarity of $B^H$ and the formula $\E[(Y_1^2-1)(Y_2^2-1)]=2\E[Y_1Y_2]^2$ for jointly normal vector $(Y_1,Y_2)$
with standard normal marginal distribution, we conclude that 
\[
\lambda_{v_1,v_2}^k = 2\left(1+ \sum_{j\in \mathbb Z\setminus \{0\}} \rho_k^{v_1,v_2} (j)^2\right),
\]
where the latter series is finite for all $k\geq 2$ and also for $k=1$ if $H=\alpha + 1/2<3/4$ holds. The condition
$H<3/4$ is well known in the framework of Breuer-Major central limit theorems for quadratic functionals (see \cite{BM83}).
This condition directly translates to $\alpha <1/4$. However, we require an additional restriction $\alpha <0$
when $k=1$ in Theorem \ref{th2} due to a certain bias, which might affect the central limit theorem. \qed
\end{rem}

\begin{rem} \label{rem5} \rm
Theorem \ref{th2} deals with realised quadratic variation only, since it is sufficient for the estimation of the smoothness
parameter $\alpha$ as we will see below. However, we do think that the asymptotic theory can be extended to functionals of the 
type
\[
V(X,h,k,v\Delta_n)_t := \sum_{i=vk}^{[t/\Delta_n]} h\left ( \frac{\Delta_{i,k}^{n,v} X}{\tau_{k}(v\Delta_n)} \right),
\]
where $h\in C^1(\R)$ is an even function. The main step of the proof is the approximation 
\[
\Delta_{i,k}^{n,v} X \approx \sum_{j=0}^l \sigma_{i\Delta_n-\theta_j} \Delta_{i,k}^{n,v} G^{(j)},
\]
where $\Delta_{i,k}^{n,v} G^{(j)}$, $j=0,\ldots,l$, are certain Gaussian random variables. 
Using Bernstein's blocking technique, which amounts in freezing
the intermittency process $\sigma$ in the beginning of sub-blocks, the asymptotic behaviour of the statistic $V(X,h,k,v\Delta_n)$
is determined by the functional 
\[
Q(z,\widetilde{h},k,v\Delta_n )_t := 
\sum_{i=vk}^{[t/\Delta_n]} \widetilde{h}\left ( z_0 \frac{\Delta_{i,k}^{n,v} G^{(0)}}{\tau_{k}(v\Delta_n)},
\ldots,  z_l \frac{\Delta_{i,k}^{n,v} G^{(l)}}{\tau_{k}(v\Delta_n)} \right), \qquad z\in \R^{l+1},
\]  
where $\widetilde{h}\in C^1(\R^{l+1})$. The central limit theorem for a standardized
version of $Q(\cdot,\widetilde{h},k,v\Delta_n )$ relies
on the stable convergence of finite dimensional distributions and tightness. The convergence of finite dimensional distributions
is a classical setting of Breuer-Major central limit theorem. It can be shown via method of moments or using more
modern methods of Malliavin calculus (see \cite{NP05,PT05} among others). We remark that in the case $h(x)=x^p$, where $p$
is an even number, we do not need to consider the process $Q(\cdot,\widetilde{h},k,v\Delta_n )$ and the proof becomes simpler
due to binomial formula. \qed  
\end{rem}

\section{The ratio statistic} \label{sec4}
\setcounter{equation}{0}
\renewcommand{\theequation}{\thesection.\arabic{equation}}
The smoothness parameter $\alpha$ defined at \eqref{min} describes the H\"older continuity index of $X$,
i.e. $X$ is H\"older continuous of any order smaller than $H=\alpha +1/2$. In the context of turbulence 
modeling the parameter $\alpha$ is connected to the so called Kolmogorov's $2/3$-law (see \cite{K41}).
It predicts that $\alpha \approx -1/6$ (or, in other words, $2(\alpha +1/2)\approx 2/3$). From this perspective
it is important to construct a consistent estimator of $\alpha$ to check if $\mathcal{BSS}$ models
adequately describe the physical laws.  

The next lemma is  crucial for estimating $\alpha$. 

\begin{lem} \label{lem1}
Assume that conditions (A) and \eqref{robustness} hold. When $k=1$ we further assume that $\alpha_j \in (-1/2,0)$
for all $0\leq j\leq l$. 
Then we obtain
\begin{align} \label{tauformula}
\tau_{k}(v\Delta_n)^2 = (v\Delta_n)^{2\alpha +1} \sum_{j=1}^l \|h_j\|_{\mathbb{L}^2(\R)}^2
1_{j\in \mathcal{A}} + o(\Delta_n^{2\alpha +3/2}), 
\end{align} 
where the functions $h_j$ were defined in \eqref{hdef}.
\end{lem}
Now, Lemma \ref{lem1} and Theorem \ref{th1} provide a 
direct way of estimating the scaling parameter $\alpha$. Indeed, we observe that 
\begin{align*}
S_n:=\frac{QV(X,k,2\Delta_n)_t}{QV(X,k,\Delta_n)_t} \toop 2^{2\alpha +1},
\end{align*}
for any fixed $t>0$. Thus,  a consistent estimator of $\alpha$ is given via
\begin{align} \label{hatalpha}
\widehat{\alpha}_n=\frac 12 \left( \log_2 \left( \frac{QV(X,k,2\Delta_n)_t}{QV(X,k,\Delta_n)_t}  \right)  -1
\right) \toop \alpha ,
\end{align}
where $\log_2$ denotes the logarithm at basis $2$. We remark that this is exactly the same estimator as proposed
in \cite{BCP11,BCP13} for $\mathcal{BSS}$ processes with a single singularity at $0$. A feasible central limit theorem 
for $\widehat{\alpha}_n$  is obtained as follows. Note that the result below requires the knowledge of
the singularity points $\theta_j$.

\begin{theo} \label{th3}
Assume that conditions of Theorem \ref{th2} hold and choose $t<\min_{1\leq j\leq l}(\theta_j-\theta_{j-1})$. \\ \\
(i) Define
\begin{align} \label{qq}
QQ(X,k,v\Delta_n)_t := 
\sum_{i=vk}^{[t/\Delta_n]} (\Delta_{i,k}^{n,v} X)^4.
\end{align}
Then we obtain that
\[
\frac{\Delta_n}{\tau_{k}(v\Delta_n)^4} 
QQ(X,k,v\Delta_n)_t \ucp 3 \int_0^t \left( \int_0^{\infty} \sigma_{s-\theta}^2 ~\pi_k(d\theta) \right)^2 ds
\]
(ii) Furthermore,  we have for any fixed $t>0$
\begin{align} \label{feasibleclt}
\frac{2\log(2) QV(X,k,\Delta_n)_t  (\widehat{\alpha}_n - \alpha)}
{\sqrt{\frac 13 QQ(X,k,\Delta_n)_t (-1,1) \Lambda_k^n (-1,1)^{\star}}} \schw \mathcal N(0,1),
\end{align}
where $\log$ denotes the logarithm at basis $e$, $x^{\star}$ is the transpose of $x$ and the matrix $\Lambda_k^n$
is defined as $\Lambda_k$ in \eqref{lambda}, where the unknown parameter $\alpha$ is replaced by $\widehat{\alpha}_n$
(recall that due to Remark \ref{rem4} the matrix $\Lambda_k$ is a function of $\alpha$).
\end{theo} 
\textit{Proof.} Here we demonstrate the proof of part (ii), while part (i) will be proved in Section \ref{sec5}.
First of all, we remark that under the condition $t<\min_{1\leq j\leq l}(\theta_j-\theta_{j-1})$ we may apply
the result of Theorem \ref{th2}. 
Define the functions $r(x,y)=2^{2\alpha +1} \frac yx$ and $h(x)=\frac 12 (\log_2(x) -1)$. Lemma \ref{lem1} implies
that 
\[
\frac{\tau_{k}(2\Delta_n)^2}{\tau_{k}(\Delta_n)^2} = 2^{2\alpha +1} +o(\Delta_n^{1/2}).
\]
Hence,
\[
S_n = r\left(\frac{\Delta_n}{\tau_{k}(\Delta_n)^2} QV(X,k,\Delta_n)_t, \frac{\Delta_n}{\tau_{k}(2\Delta_n)^2} QV(X,k,2\Delta_n)_t
\right)  +o_{\mathbb P}(\Delta_n^{1/2}).
\]
Putting things together we conclude that 
\[
\Delta_n^{-1/2} (\widehat{\alpha}_n - \alpha) = h\circ r 
\left(\frac{\Delta_n}{\tau_{k}(\Delta_n)^2} QV(X,k,\Delta_n)_t, \frac{\Delta_n}{\tau_{k}(2\Delta_n)^2} QV(X,k,2\Delta_n)_t
\right)  +o_{\mathbb P}(1).
\]
Applying Theorem \ref{th2} and delta method for stable convergence we deduce that 
\[
\Delta_n^{-1/2} (\widehat{\alpha}_n - \alpha) \stab \mathcal{MN}(0, V^2),
\]
where $\mathcal{MN}(0, V^2)$ denotes a mixed normal variable with mean $0$ and conditional variance $V^2$ defined by
\[
V^2:= \frac{(-1,1) \Lambda_k (-1,1)^{\star} 
\int_0^t \left( \int_0^{\infty} \sigma_{s-\theta}^2 ~\pi_k(d\theta) \right)^2 ds}{\Big(2\log (2) QV(X,k)_t \Big)^2 },
\]
where the matrix $\Lambda_k$ is defined by \eqref{lambda}. Notice that $\Lambda_k$ is a continuous function 
of $\alpha$ due to Remark \ref{rem4}. Hence,
\[
\Lambda_k^n \toop \Lambda_k. 
\]
The two other random quantities involved in the definition of $V^2$ can be directly estimated via part (i) of Theorem
\ref{th3} and Theorem \ref{th1}. Consequently, the properties of stable convergence imply part (ii) of Theorem
\ref{th3}. \qed 
\\ \\
Note that the standardized statistic in \eqref{feasibleclt} is feasible as it does not require the knowledge of the weight 
function $g$. We remark that \eqref{feasibleclt} coincides with the statistic presented in \cite[Proposition 4.2]{CHPP13}
in the framework of a single singularity at $0$. This demonstrates that the test statistic in \eqref{feasibleclt}
is robust to model misspecification within the setting of assumption (A) and condition \eqref{robustness}. 
In the context of turbulence modeling
this is a very important property.

\section{Proofs} \label{sec5}
\setcounter{equation}{0}
\renewcommand{\theequation}{\thesection.\arabic{equation}}

\subsection{Proof of Proposition \ref{prop1} and Lemma \ref{lem1}} \label{sec5.1}
We first prove Lemma \ref{lem1} as its proof essentially implies Proposition \ref{prop1}. Throughout this section all
positive constants are denoted by $C$, or $C_p$ if they depend on an external parameter $p$, although they may change
from line to line. \\ \\
\textit{Proof of Lemma \ref{lem1}.} We assume without loss of generality that $l=2$, $\alpha_0=\alpha_1=\alpha$
and $\alpha_2-\alpha >1/4$ (since condition \eqref{robustness} was assumed). Moreover, let $v=1$. Recall the
identity
\[
\tau_{k}(\Delta_n)^2 = \|\Delta_{k}^{n,1} g\|_{\mathbb{L}^2(\R)}^2
\]
(cf. Section \ref{sec3.1}). We consider the decomposition
\begin{align}
\|\Delta_{k}^{n,1} g\|_{\mathbb{L}^2(\R)}^2 &= \sum_{j=0}^2 \int_{\theta_j-\delta}^{\theta_j+\delta} \Delta_{k}^{n,1} g(x)^2 dx
+ \int_{\delta}^{\theta_1-\delta} \Delta_{k}^{n,1} g(x)^2 dx+ 
\int_{\theta_1+\delta}^{\theta_2-\delta} \Delta_{k}^{n,1} g(x)^2 dx  \nonumber \\[1.5 ex]
\label{taudec}&+ \int_{\theta_2+\delta}^{\infty} \Delta_{k}^{n,1} g(x)^2 dx.
\end{align}
We will now show that 
\begin{align} \label{domterms}
\int_{\theta_j-\delta}^{\theta_j+\delta} \Delta_{k}^{n,1} g(x)^2 dx 
= \Delta_n^{2\alpha +1}  \|h_j\|_{\mathbb{L}^2(\R)}
+ o(\Delta_n^{2\alpha +3/2}), \qquad j=0,1,
\end{align}
and all other terms in the decomposition are $o(\Delta_n^{2\alpha +3/2})$ under the assumptions of Lemma \ref{lem1}.
We start with the negligibility of the three last terms in \eqref{taudec}. Using Taylor expansion of order $k$ and integrability
condition (A)(iv), we immediately conclude that  
\begin{align*}
\int_{\delta}^{\theta_1-\delta} \Delta_{k}^{n,1} g(x)^2 dx &= O(\Delta_n^{2k}), \quad
\int_{\theta_1+\delta}^{\theta_2-\delta} \Delta_{k}^{n,1} g(x)^2 dx = O(\Delta_n^{2k}), \\
\int_{\theta_2+\delta}^{\infty} \Delta_{k}^{n,1} g(x)^2 dx &= O(\Delta_n^{2k}),
\end{align*}
so all these terms are $o(\Delta_n^{2\alpha +3/2})$ under assumptions of Lemma \ref{lem1}.
Now, we show \eqref{domterms} for $j=0$; the case $j=1$ works similarly. Proving  this statement for $j=0$ 
essentially means that we can replace $f_0(x)$ involved in the integral by the constant $f_0(0)$. Let $\varepsilon >0$
be small enough with $\varepsilon >> \Delta_n$. Using again 
Taylor expansion of order $k$ and integrability condition (A)(iv), we conclude that (recall that $g(x)=0$
for $x\leq 0$)
\[
\int_{-\delta}^{\delta} \Delta_{k}^{n,1} g(x)^2 dx = 
\int_{0}^{\varepsilon} \Delta_{k}^{n,1} g(x)^2 dx + O(\Delta_n^{2k}\varepsilon^{2(\alpha -k)+1}).  
\] 
When we replace the function $f_0$ that appears in the latter integral by a constant $f_0(0)$, we deduce by substitution
$x=\Delta_n y$
\begin{align*}
f_0(0)^2 \int_{0}^{\varepsilon} (\Delta_{k}^{n,1} (x^{\alpha}))^2 dx &= \Delta_n^{2\alpha +1}
\int_{0}^{\varepsilon /\Delta_n} h_0(y)^2 dy \\[1.5 ex]
&= \Delta_n^{2\alpha +1}
\int_{0}^{\infty } h_0(y)^2 dy + O(\Delta_n^{2k}\varepsilon^{2(\alpha -k)+1}), 
\end{align*}
since $|h_0(x)|^2\leq C |x|^{2(\alpha-k)}$ for large $x$ and $2(\alpha-k)<-1$. Note that the dominating term 
is exactly the one given in \eqref{domterms}. Now, let us evaluate the difference
\begin{align*}
\int_{0}^{\varepsilon} \Delta_{k}^{n,1} g(x)^2 dx - f_0(0)^2 \int_{0}^{\varepsilon} (\Delta_{k}^{n,1} (x^{\alpha}))^2 dx &=
\int_{0}^{k\Delta_n} \Delta_{k}^{n,1} g(x)^2 - f_0(0)^2 (\Delta_{k}^{n,1} (x^{\alpha}))^2 dx
\\[1.5 ex]
&+ \int_{k\Delta_n}^{\varepsilon} \Delta_{k}^{n,1} g(x)^2 - f_0(0)^2 (\Delta_{k}^{n,1} (x^{\alpha}))^2 dx.
\end{align*} 
Using differentiability of $f_0$ we immediately conclude that 
\[
\int_{0}^{k\Delta_n} |\Delta_{k}^{n,1} g(x)^2 - f_0(0)^2 (\Delta_{k}^{n,1} (x^{\alpha}))^2| dx = O(\Delta_n^{2\alpha +2}).
\]
The other integral has to be treated differently. In the following we present the computations only for $k=1,2$
(in fact, the case $k\geq 3$ is easier to treat). We start with $k=1$. Using binomial rule, differentiability of $f_0$
and substitution, we conclude that
\[
\Big|\int_{\Delta_n}^{\varepsilon} \Delta_{1}^{n,1} g(x)^2 - f_0(0)^2 (\Delta_{1}^{n,1} (x^{\alpha}))^2 dx \Big|\leq 
C \Delta_n^{2\alpha +2} \int_{1}^{\varepsilon /\Delta_n} |h_0(y)| y^{\alpha+1} dy 
\]
Using the inequality $|h_0(x)|\leq C |x|^{\alpha-1}$ for large $x$, we deduce that
\[
\Delta_n^{2\alpha +2} \int_{k}^{\varepsilon /\Delta_n} |h_0(y)| y^{\alpha+1} dy \leq C
\Delta_n \varepsilon^{2\alpha +1}
\]
Setting $\varepsilon =\Delta_n^{1/2}$, we deduce that all involved small order terms are $o(\Delta_n^{2\alpha +3/2})$
when $\alpha <0$. Now, we consider the case $k=2$. Since $f_0$ is twice continuously differentiable, we may apply the Taylor
expansion to
\[
f_0(x+m\Delta_n) = f_0(0) + mf'_0(0)\Delta_n + \frac 12 m^2 \Delta_n^2 f''_0(x_m),  
\]
where $x_m\in (0,x+m\Delta_n)$ and $m=0,1,2$. Using the above Taylor expansion and the binomial formula,
we deduce  that 
\[
\Big| \int_{k\Delta_n}^{\varepsilon} \Delta_{k}^{n,1} g(x)^2 - f_0(0)^2 (\Delta_{k}^{n,1} (x^{\alpha}))^2 dx \Big|
\leq C(q_1(n,\varepsilon ) + q_2(n,\varepsilon ) +q_3(n,\varepsilon )),
\]
where
\begin{align*}
q_1(n,\varepsilon )&= \int_{k\Delta_n}^{\varepsilon} h_0(x)^2 x dx, \\[1.5 ex]
q_2(n,\varepsilon )&= \Delta_n\int_{k\Delta_n}^{\varepsilon} |h_0(x)| |(x+2\Delta_n)^\alpha - (x+\Delta_n)^\alpha| dx, \\[1.5 ex]
q_3(n,\varepsilon )&= \int_{k\Delta_n}^{\varepsilon} |h_0(x)| x^{\alpha +2} dx.
\end{align*}
Applying the substitution $x=\Delta_n y$, we get
\[
q_1(n,\varepsilon ) = O(\Delta_n^{2\alpha +2}), \qquad q_2(n,\varepsilon ) = O(\Delta_n^{2\alpha +2}), \qquad 
q_3(n,\varepsilon ) = O(\Delta_n^{2} \varepsilon^{2\alpha +1}).
\]
Setting now $\varepsilon = \Delta_n^{1/2}$, we conclude that all second order terms are $o(\Delta_n^{2\alpha +3/2})$. 

Finally, let us treat the case $j=2$. Since 
\[
\int_{\theta_2-\delta}^{\theta_2+\delta} \Delta_{k}^{n,1} g(x)^2 dx =O (\Delta_n^{2\alpha_2+1}) 
\]
as shown above and $\alpha_2-\alpha >1/4$, we see that this term is $o(\Delta_n^{2\alpha + 3/2})$.
Consequently, we obtain the assertion of Lemma \ref{lem1}. \qed \\ \\
\textit{Proof of Proposition \ref{prop1}.} The assertion of Proposition \ref{prop1} now easily follows from the proof
of Lemma \ref{lem1}. First of all, it implies that
\[
\|\Delta_{k}^{n,v} g\|_{\mathbb{L}^2(\R)}^2 = (v\Delta_n)^{2\alpha +1} 
\sum_{j=0}^l \|h_j\|_{\mathbb{L}^2(\R)}^2
1_{j\in \mathcal{A}} + o(\Delta_n^{2\alpha +1}),
\] 
even without condition \eqref{robustness}. On the other hand the proof of Lemma \ref{lem1} also implies that
\[
\int_{\theta_j-\varepsilon}^{\theta_j+\varepsilon} \Delta_{k}^{n,v} g(x)^2 dx = (v\Delta_n)^{2\alpha +1}
\|h_j\|_{\mathbb{L}^2(\R)}^2 + o(\Delta_n^{2\alpha +1}) \qquad \text{if } j\in \mathcal{A},
\]
for any $\varepsilon < \min_{1\leq i\leq l}(\theta_i-\theta_{i-1})$, and 
\[
\int_{a}^{b} \Delta_{k}^{n,v} g(x)^2 dx =  o(\Delta_n^{2\alpha +1})
\]
if the interval $[a,b]$ does not contain any $\theta_j$ with $j\in \mathcal{A}$. This completes the proof of Proposition
\ref{prop1}. \qed

\subsection{Some preliminaries} \label{sec5.2}
Before we prove the main results of the paper, we start with some preliminaries. 
We remark that the intermittency process $\sigma$ is c\'adl\'ag, hence $\sigma_{-}$ is locally bounded.
Since our Theorems \ref{th1}, \ref{th2} and \ref{th3} are stable under localization (cf. \cite{BGJPS06}),
we may and will assume that $\sigma$ is bounded on compact intervals.

Recalling the notation of \eqref{filterg} we introduce the following Gaussian random variables

\begin{align} \label{gj}
\Delta_{i,k}^{n,v} G^{(j)}:= \int_{i\Delta_n-\theta_j-\delta}^{i\Delta_n-\theta_j+\delta} \Delta_{k}^{n,v} g(i\Delta_n-s) W(ds),
\qquad j=0,\ldots,l,
\end{align}
where the constant $\delta >0$ was defined in (A). Notice that the above Gaussian variables
are independent for different $j$'s when computed at the same stage $i\Delta_n$. 
In the first step of our proofs we will show the approximation 
\begin{align} \label{xapprox}
\Delta_{i,k}^{n,v} X \approx \sum_{j=0}^l \sigma_{(i-vk)\Delta_n-\theta_j} \Delta_{i,k}^{n,v} G^{(j)}.
\end{align}
The ideas behind the proofs of Theorems \ref{th1}, \ref{th2} and \ref{th3} follow a similar structure as presented
in \cite{BCP11,BCP13}, although the situation is more complex due to multiple singularities of the weight function $g$.
First of all, we will use a blocking technique, which amounts in considering a subdivision of the interval $[0,t]$ into
equidistant sub-blocks and freezing the intermittency process $\sigma$ within each sub-block. In a second step, 
we will prove joint limit theorems over the sub-blocks applying Malliavin calculus and properties of stable convergence.  

We start with the limit theory for the Gaussian variables $\Delta_{i,k}^{n,v} G^{(j)}$, which has been essentially treated in
\cite{BCP11}. Define
\begin{align} \label{rnj}
\tau_{k,j}(v\Delta_n)^2:= \E[(\Delta_{i,k}^{n,v} G^{(j)})^2], \qquad r_{k,n}^{v_1,v_2,j_1,j_2}(q):=\text{corr}
\Big(\Delta_{1,k}^{n,v_1} 
G^{(j_1)},
\Delta_{1+q,k}^{n,v_2} G^{(j_2)} \Big).
\end{align}
We consider the statistics
\begin{align} \label{qvg}
QV(k,v\Delta_n)_t^{j_1,j_2} &:= \Delta_n \sum_{i=vk}^{[t/\Delta_n]} \frac{\Delta_{i,k}^{n,v} G^{(j_1)} \Delta_{i,k}^{n,v} G^{(j_2)}}
{\tau_{k,j_1}(v\Delta_n) \tau_{k,j_2}(v\Delta_n)}, \qquad 0\leq j_1,j_2\leq l,
\end{align}
and set
\begin{align} \label{rhov}
\rho^{v_1,v_2,j}_{k}(q):= \text{corr}(\Delta_{1,k}^{v_1} B^{H_j}, \Delta_{1+q,k}^{v_2} B^{H_j}), 
\end{align}
where $B^{H_j}$ is a fractional Brownian motion with Hurst parameter $H_j=\alpha_j+1/2$ and the quantity 
$\Delta_{i,k}^v B^{H_j}$ is defined in \eqref{fbmnoise}.
The next result is essentially a combination of \cite[Theorems 1 and 2]{BCP11} and \cite[Section 2]{BCP13}.

\begin{theo} \label{th5}
Assume that condition (A) holds. \\ \\
(i) We have the convergence 
\begin{align} \label{qvglln}
QV(k,v\Delta_n)_t^{j_1,j_2} \ucp QV(k)_t^{j_1,j_2}:=\delta_{j_1,j_2}t, \qquad j_1,j_2\in \mathcal A,
\end{align}
where $\delta_{j_1,j_2}=1$ when $j_1=j_2$ and $0$ otherwise. \\ \\
(ii) When $k=1$ we further assume that $\alpha_j <0$ for all $0\leq j\leq l$. Then 
\begin{align} \label{qvgclt}
\Delta_n^{-1/2} \left( QV(k,v\Delta_n)^{j_1,j_2} - QV(k)^{j_1,j_2} \right)_{j_1,j_2\in \mathcal A,j_1\leq j_2}^{v=1,2} \stab
V=(V_{k,v}^{j_1,j_2})_{j_1,j_2\in \mathcal A,j_1\leq j_2}^{v=1,2},
\end{align}
on $\mathbb D^{|\mathcal{A}|(|\mathcal{A}|+1)}([0, \min_{1\leq j\leq l}(\theta_j-\theta_{j-1})])$,
where $V$ is an $\mathcal F$-conditional Gaussian martingale, defined on an extension $(\Omega', \mathcal
F', \mathbb P')$ of the original probability space. The   $\mathcal F$-conditional  covariance 
structure is given as
\begin{align*}
\E'[V_{k,v}^{j_1,j_2}(t) V_{k,v'}^{j'_1,j'_2}(s)]&= 0 \qquad \text{when } (j_1,j_2) \not =(j'_1,j'_2), \\
\E'[V_{k,v}^{j,j}(t) V_{k,v'}^{j,j}(s)]&= 2\min\{t,s\} \left(1+ \sum_{q\in \mathbb Z\setminus \{0\}} \rho^{v,v',j}_{k}(q)^2\right), \\
\E'[V_{k,v}^{j_1,j_2}(t) V_{k,v'}^{j_1,j_2}(s)]&= \min\{t,s\} 
\left(1+ \sum_{q\in \mathbb Z\setminus \{0\}} \rho^{v,v',j_1}_{k}(q) \rho^{v,v',j_2}_{k}(q)\right) \quad \text{for }j_1 \not =j_2.
\end{align*} 
\end{theo}
We remark that for $j,j_1,j_2\in \mathcal A$ we immediately conclude that 
\begin{align} \label{lambdaid}
\E'[V_{k,v}^{j,j}(1) V_{k,v'}^{j,j}(1)] &= \lambda^k_{v,v'},  \\
\E'[V_{k,v}^{j_1,j_2}(1) V_{k,v'}^{j_1,j_2}(1)]&= \frac 12 \lambda^k_{v,v'} \quad \text{for }j_1 \not =j_2, \nonumber
\end{align}
where $\lambda^k_{v,v'}$ is defined by \eqref{lambda}. \\ \\
\textit{Proof.} We divide the proof into several steps. \\ \\
\textit{Step 1.} We will not work with the random variables $\Delta_{i,k}^{n,v} G^{(j)}$ directly, but with their approximations. 
For any $0\leq j\leq l$, we define a new Gaussian process
\begin{align*}
\widetilde{G}_t^{(j)}:= \int_{\R} \widetilde g^{(j)} (t-s) W(ds),
\end{align*} 
where the function $\widetilde g^{(j)}$ ($j\geq 1$) is given via
\begin{align*}
\widetilde g^{(j)} (x) = \widetilde f_j(x) |x- \theta_j|^{\alpha_j},  
\end{align*} 
with $\widetilde f_j=f_j$ on $x\in (\theta_j - \delta/2, \theta_j + \delta/2)$, 
$\widetilde f_j = 0$ outside of the interval $(\theta_j - \delta, \theta_j + \delta)$, and 
$\widetilde f_j \in C^k(\R )$ (the function $\widetilde f_{0}$ is defined similarly).
Then, if we consider the $k$-th order increments $\Delta_{i,k}^{n,v} \widetilde G^{(j)}$  of $\widetilde{G}$ at frequency
$v\Delta_n$, we readily deduce that 
\begin{align} \label{apperror}
\E[(\Delta_{i,k}^{n,v} \widetilde G^{(j)} - \Delta_{i,k}^{n,v} G^{(j)})^2]\leq C \Delta_n^{2k} 
\end{align}   
due to assumption (A). Now, let us define the statistics
\begin{align} \label{tildeqvg}
\widetilde{QV}(k,v\Delta_n)_t^{j_1,j_2} &:= \Delta_n \sum_{i=vk}^{[t/\Delta_n]} \frac{\Delta_{i,k}^{n,v} \widetilde G^{(j_1)} 
\Delta_{i,k}^{n,v} \widetilde G^{(j_2)}}
{\widetilde{\tau}_{k,j_1}(v\Delta_n) \widetilde{\tau}_{k,j_2}(v\Delta_n)}, \qquad 
\widetilde{\tau}_{k,j}(v\Delta_n)^2:= \E[(\Delta_{i,k}^{n,v} \widetilde G^{(j)})^2],
\end{align}
for $0\leq j_1,j_2\leq l$. Then, due to \eqref{apperror}, $\widetilde{\tau}_{k,j}(v\Delta_n)/ \tau_{k,j}(v\Delta_n)\rightarrow 1$
and 
\begin{align} \label{pnnegl} 
\widetilde{QV}(k,v\Delta_n)_t^{j_1,j_2} - QV(k,v\Delta_n)_t^{j_1,j_2} \ucp 0,
\end{align}
and also 
\begin{align} \label{cltnegl} 
\Delta_n^{-1/2} \left( \widetilde{QV}(k,v\Delta_n)_t^{j_1,j_2} - QV(k,v\Delta_n)_t^{j_1,j_2} \right) \ucp 0
\end{align}
under the assumption of Theorem \ref{th5}(ii), which is due to Cauchy-Schwarz inequality. Thus, it suffices to prove the asymptotic
theory for the statistics $\widetilde{QV}(k,v\Delta_n)_t^{j_1,j_2}$. \qed \\ \\
\textit{Step 2.} In this step we analyze the correlation structure of the increments $\Delta_{i,k}^{n,v} \widetilde G^{(j)}$.  
We define
\begin{align} \label{tilde rnj}
\widetilde{r}_{k,n}^{v_1,v_2,j_1,j_2}(q):=\text{corr}
\Big(\Delta_{1,k}^{n,v_1} 
\widetilde G^{(j_1)},
\Delta_{1+q,k}^{n,v_2} \widetilde G^{(j_2)} \Big).
\end{align}
The next proposition describes the asymptotic behaviour of the correlation function $\widetilde r_{k,n}^{v_1,v_2,j_1,j_2}(q)$.
\begin{prop}
Assume that condition (A) holds. Then we obtain that
\begin{align} \label{conv1}
\widetilde r_{k,n}^{v_1,v_2,j,j}(q) \rightarrow 
\rho^{v_1,v_2,j}_{k}(q), 
\end{align}
where $\rho^{v_1,v_2,j}_{k}(q)$ is defined at \eqref{rhov}. Furthermore, for any $\epsilon >0$ there exists $C>0$
such that
\begin{align} \label{bound1}
|\widetilde r_{k,n}^{v_1,v_2,j,j}(q)|\leq C |q|^{2H_j-2k-\epsilon}.
\end{align}
Likewise, if $j_1,j_2\in \mathcal A$, $j_1>j_2$ then for any $\epsilon >0$ there exists $C>0$
such that 
\begin{align} \label{bound2}
|\widetilde r_{k,n}^{v_1,v_2,j_1,j_2}(q)|\leq C |q+\Delta_n^{-1}(\theta_{j_1}-\theta_{j_2})|^{2\alpha +1-2k-\epsilon}.
\end{align} 
\end{prop}
\textit{Proof.} We start with the proof of \eqref{conv1} and \eqref{bound1}. Without loss of generality we prove it only for the case
$k=1$, $v_1=v_2=1$ and $j=0$; the rest follows by similar arguments. For simplicity we set $ \widetilde r_{n}(q):=
\widetilde r_{1,n}^{1,1,0,0}(q)$, $\rho(q):= \rho^{1,1,0}_{1}(q)$ and $\Delta_{i}^{n} \widetilde G^{(0)}:=
\Delta_{i,1}^{n,1} \widetilde G^{(0)}$. Observe
that
\begin{align*}
\text{cov} \left(\Delta_{1}^{n} \widetilde G^{(0)}, \Delta_{1+q}^{n} \widetilde G^{(0)} \right) &= 
\int_{(1+q)\Delta_n-\delta}^{\Delta_n+\delta} \{\widetilde g^{(0)}(\Delta_n-s)- \widetilde g^{(0)}(-s)\} \\[1.5 ex]
&\times \{\widetilde g^{(0)}((1+q)\Delta_n-s)- \widetilde g^{(0)}(q\Delta_n-s)\} ds.
\end{align*} 
Now, recalling that $\widetilde g^{(0)}(x)= x^{\alpha_0} f_0(x)$ for $x\in (0,\delta/2)$, we conclude as in the proof of 
Lemma \ref{lem1}
\begin{align*}
\Delta_n^{-(2\alpha_0 +1)}\text{cov} \left(\Delta_{1}^{n} \widetilde G^{(0)}, \Delta_{1+q}^{n} \widetilde G^{(0)} \right) &\rightarrow 
f_0(0)^2
\int_{\R} \{ (1-s)_{+}^{\alpha_0}- (-s)_{+}^{\alpha_0} \} \\[1.5 ex]
&\times \{ (q+1-s)_{+}^{\alpha_0}- (q-s)_{+}^{\alpha_0} \} ds.
\end{align*}
The latter limit is, up to a factor $f_0(0)^2$, the covariance function of a (non-standard) fractional Brownian noise
$(\widetilde{B}_i^{\alpha_0+1/2} - \widetilde{B}_{i-1}^{\alpha_0+1/2} )_{i\geq 1}$, where $\widetilde{B}^{\alpha_0+1/2}$ 
is defined as
\[
\widetilde{B}_t^{\alpha_0+1/2} := \int_{\R} \{(t-s)_{+}^{\alpha_0} -(-s) _{+}^{\alpha_0}\} W(ds).
\]
Thus, using again Lemma \ref{lem1}, we deduce that
\[
\widetilde r_{n}(q) \rightarrow \rho (q) \qquad \text {as } n\rightarrow \infty, 
\] 
which completes the proof of \eqref{conv1}.

Now, we define the function
\[
\widetilde R_u:= \int_{-\delta}^{\delta} (\widetilde g^{(0)}(-s) - \widetilde g^{(0)}(-u-s))^2  ds
\]  
and note that $\widetilde R_{\Delta_n}= \widetilde{\tau}_{1,0}(\Delta_n)^2$. According 
to \cite[Lemma 1]{BCP09} and conditions (A1)-(A3) therein,
it is sufficient to show that 
\[
\widetilde R_u = u^{2\alpha_0+1} Z_u, \qquad u> 0,
\] 
where $Z\in C^2(0,\infty )$ and $\lim_{u\rightarrow 0} Z_u \not =0$, to conclude  \eqref{bound1}. Observe that
for $u<\delta$ 
\begin{align*}
\widetilde R_u &= 
\int_{0}^{u} s^{2\alpha_0} \widetilde f^2_0(s)  ds +  \int_{u}^{\delta} (s^{\alpha_0} \widetilde f_0(s) -
(s-u)^{\alpha_0} \widetilde f^2_0(s-u))^2 ds 
\\[1.5 ex]
&= u^{2\alpha_0+1}  \left( \int_{0}^{1} x^{2\alpha_0} \widetilde f^2_0(ux) dx + \int_{1}^{\delta/u} 
(x^{\alpha_0} \widetilde f_0(ux) -(x-1)^{\alpha_0} \widetilde f_0^2(ux-u))^2 dx \right)  \\[1.5 ex]
&= u^{2\alpha_0+1} Z_u. 
\end{align*}
Now, $Z\in C^2(0,\infty )$ since $\widetilde f_0\in C^2(0,\delta)$; see condition (A)(iii). Furthermore,
\[
\lim_{u\rightarrow 0} Z_u= f_0^2(0) \left( \int_{0}^{1} x^{2\alpha_0}  dx + \int_{1}^{\infty} 
(x^{\alpha_0}  -(x-1)^{\alpha_0} )^2 dx \right),
\]
where the limit is finite since $\alpha_0<1/2$ and strictly positive because $f_0(0)\not=0$ (see again (A)(iii)).
Thus,  \eqref{bound1} follows.

Finally, we need to prove \eqref{bound2}. We assume without loss of generality that $k=1$, $v_1=v_2=1$ and $\theta_{j_1}>
\theta_{j_2} > 0$. For $q<0$ such that $(1+q)\Delta_n-\theta_{j_2}-\delta \in 
[\Delta_n-\theta_{j_1}-\delta, \Delta_n-\theta_{j_1}-\delta)$, we obtain that 
\begin{align*}
&\text{cov} \left(\Delta_{1}^{n} \widetilde G^{(j_1)}, \Delta_{1+q}^{n} \widetilde G^{(j_2)} \right) \\[1.5 ex]
&= \int_{(1+q)\Delta_n-\theta_{j_2}-\delta}^{\Delta_n-\theta_{j_1}+\delta} \{ \widetilde g^{(j_1)}(\Delta_n-s)- 
\widetilde g^{(j_1)}(-s)\}
\{\widetilde g^{(j_2)} ((1+q)\Delta_n-s)- \widetilde g^{(j_2)}(q\Delta_n-s)\} ds \\[1.5 ex]
&= \int_{(1+q)\Delta_n-\theta_{j_2}-\delta}^{\Delta_n-\theta_{j_1}+\delta} \{\widetilde f_{j_1}(\Delta_n-s) |\Delta_n-s-\theta_{j_1}|^\alpha 
- \widetilde f_{j_1}(-s) |s+\theta_{j_1}|^\alpha \} \\[1.5 ex]
& \times \{ \widetilde f_{j_2}((1+q)\Delta_n-s) |(1+q)\Delta_n-s-\theta_{j_2}|^\alpha 
- \widetilde f_{j_2}(q\Delta_n-s) |q\Delta_n-s-\theta_{j_2}|^\alpha \} ds,
\end{align*}
where we recall that $\alpha_{j_1}=\alpha_{j_2}=\alpha$. In the next step we compare this expression with
the following covariance
\begin{align*}
&\text{cov} \left(\Delta_{1}^{n} \widetilde G^{(j_1)}, \Delta_{1+\bar{q}}^{n} \widetilde G^{(j_1)} \right) \\[1.5 ex]
&= \int_{(1+\bar{q})\Delta_n-\theta_{j_1}-\delta}^{\Delta_n-\theta_{j_1}+\delta} \{\widetilde f_{j_1}(\Delta_n-s) 
|\Delta_n-s-\theta_{j_1}|^\alpha 
- \widetilde f_{j_1}(-s) |s+\theta_{j_1}|^\alpha \} \\[1.5 ex]
& \times \{ \widetilde f_{j_1}((1+\bar{q})\Delta_n-s) |(1+\bar{q})\Delta_n-s-\theta_{j_1}|^\alpha 
- \widetilde f_{j_1}(\bar{q}\Delta_n-s) |\bar{q} \Delta_n-s-\theta_{j_1}|^\alpha \} ds.
\end{align*}
Now, by setting $\bar{q}=[q + \Delta_n^{-1}(\theta_{j_1} - \theta_{j_2})]$ and recalling that all functions $\widetilde f_j$
satisfy the same assumption (A)(iii), and keeping in mind Lemma \ref{lem1}, we conclude that 
\[
|r_{n}^{j_1,j_2}(q)|\leq C |r_{n}^{j_1,j_1}(\bar{q})|,
\] 
which implies \eqref{bound2} by applying \eqref{bound1}. \qed \\ \\
\textit{Step 3.} Due to Step 1 it suffices to prove Theorem \ref{th5} for the statistics $\widetilde{QV}(k,v\Delta_n)_t^{j_1,j_2}$. 
We start with part (i). Assertions \eqref{conv1} and \eqref{bound1} immediately imply the convergence \eqref{qvglln} by  
\cite[Theorem 1]{BCP11} (or, more precisely, by its multivariate extension).

Part (ii) essentially follows from \cite[Theorem 2]{BCP11}. First, we observe that our multivariate statistic is a functional 
of a Gaussian process. In this case it is sufficient to prove asymptotic normality for each component and to identify the covariance
structure (this is due to the results of \cite{PT05}). The asymptotic normality follows from the square summability 
of the bound in \eqref{bound1}, i.e.
\[
\sum_{q=1}^{\infty} q^{2(2\alpha +1-2k+\epsilon )} <\infty 
\]
for $\epsilon  >0$ small enough (if $k=1$ we require that $\alpha <1/4$), which is a sufficient condition for asymptotic normality
of each component due to \cite[Theorem 2]{BCP11}. Furthermore, the convergence in \eqref{conv1} easily identifies the covariance
structure of each component (see again \cite[Theorem 2]{BCP11}), hence the last two identities of Theorem \ref{th5}. 

Now, let us prove the asymptotic independence of the involved components. As before we assume without loss of generality
that $k=1$, $v_1=v_2=1$. We define
\[
\widetilde V_n^j(t)=\Delta_n^{-1/2} \left( \widetilde{QV}(1,\Delta_n)^{j,j}_t - QV(1)^{j,j}_t \right)
\] 
and show that $\E [\widetilde V_n^{j_1}(t) \widetilde V_n^{j_2}(s)] \rightarrow 0$ for $j_1,j_2\in \mathcal A$ with $j_1 \not= j_2$ 
(the asymptotic independence of all other components is shown
in exactly the same manner). Recall that $|t-s|< \min_{1\leq j\leq l}(\theta_j-\theta_{j-1})$. We deduce that 
\[
\E [\widetilde V_n^{j_1}(t) \widetilde V_n^{j_2}(s)] = 2\Delta_n \sum_{i_1=1}^{[t/\Delta_n]} \sum_{i_2=1}^{[s/\Delta_n]}
|\widetilde r_{n}^{j_1,j_2}(i_2-i_1)|^2
\]  
Now, we consider the most critical case namely  $t=s=\min_{1\leq j\leq l}(\theta_j-\theta_{j-1})=\theta_{j_1}- \theta_{j_2}$
(so $\theta_{j_1}> \theta_{j_2}$). Then the estimate \eqref{bound2} gives
\[
|\E [\widetilde V_n^{j_1}(t) \widetilde V_n^{j_2}(s)]|\leq C \Delta_n \sum_{i=-[t/\Delta_n]+1}^{[t/\Delta_n]-1} 
|i+\Delta_n^{-1}(\theta_{j_1}-\theta_{j_2})|^{2(2\alpha -1+\epsilon)} ([t/\Delta_n]-|i|) 
\]
Let $w\in (0,1)$. We conclude that
\[
\Delta_n \sum_{i=-w[t/\Delta_n]+1}^{[t/\Delta_n]-1} 
|i+\Delta_n^{-1}(\theta_{j_1}-\theta_{j_2})|^{2(2\alpha -1+\epsilon)} ([t/\Delta_n]-|i|)\leq C(1-w) \Delta_n^{-4\alpha +1 -2\epsilon}. 
\]
On the other hand we have that 
\[
\Delta_n \sum_{i=-[t/\Delta_n]+1}^{-w[t/\Delta_n]} |i+\Delta_n^{-1}(\theta_{j_1}-\theta_{j_2})|^{2(2\alpha -1+\epsilon)} 
([t/\Delta_n]-|i|)\leq C(1-w),
\] 
since $t=\theta_{j_1}- \theta_{j_2}$, $\alpha <1/4$ and $\epsilon>0$ can be chosen arbitrarily small. Hence, letting first 
$\Delta_n\rightarrow 0$ and then $w\rightarrow 1$ we obtain the desired convergence
$\E [\widetilde V_n^{j_1}(t) \widetilde V_n^{j_2}(s)] \rightarrow 0$. \qed \\ \\
Finally, let us note that due to Theorem \ref{th5}(i) we have that 
\[
\Delta_n \tau_k(\Delta_n)^{-2}\sum_{i=vk}^{[t/\Delta_n]} (\Delta_{i,k}^{n,v} G^{(j)})^2 = o_{\mathbb P} (\Delta_n^{1/2})
\qquad \forall j\not\in \mathcal{A},  
\]
due to Lemma \ref{lem1} and the condition $\alpha_j-\alpha >1/4$. On the other hand, when $j_1 \not= j_2$ and either
$j_1 \not\in \mathcal{A}$ or $j_2 \not\in \mathcal{A}$, then we conclude that  
\[
\text{var} \left( QV(k,v\Delta_n)_t^{j_1,j_2} \right) = O_{\mathbb P} (\Delta_n)
\] 
under assumptions of Theorem \ref{th5} (the arguments are similar to the proof of Theorem \ref{th5}). Thus, using again Lemma
\ref{lem1}, we conclude that
\begin{align} \label{neglig}
\Delta_n \tau_k(\Delta_n)^{-2}\sum_{i=vk}^{[t/\Delta_n]} \Delta_{i,k}^{n,v} G^{(j_1)} \Delta_{i,k}^{n,v} G^{(j_2)} = o_{\mathbb P} (\Delta_n^{1/2})
\end{align}
whenever $j_1 \not\in \mathcal{A}$ or $j_2 \not\in \mathcal{A}$, under conditions of Theorem \ref{th2}.

\subsection{Proof of Theorem \ref{th1} and Theorem \ref{th3}(i)} \label{sec5.3}
\textit{Proof of Theorem \ref{th1}.}
Below we apply a blocking technique, which means that we subdivide the interval $[0,t]$ into sub-blocks
and freeze the intermittency process within each block. We remark that the statistic $QV(X,k,v\Delta_n)_t$
is increasing in $t$ and the limiting process $QV(X,k)_t$ at \eqref{LLN} is continuous in $t$. For this reason it is sufficient
to prove pointwise convergence 
\[
\frac{\Delta_n}{\tau_{k}(v\Delta_n)^2} QV(X,k,v\Delta_n)_t \toop  QV(X,k)_t= 
\int_0^\infty \left( \int_{-\theta}^{t-\theta } \sigma_s^2 ds\right) \pi_k (d\theta)
\]
for a fixed $t>0$. 

Now, we fix a natural number $m$ and introduce the decomposition 
\begin{align} \label{decomposition}
\frac{\Delta_n}{\tau_{k}(v\Delta_n)^2} QV(X,k,v\Delta_n)_t -  QV(X,k)_t= 
A_n + B_{n,m} + C_{n,m} + D_m,
\end{align} 
where
\begin{align*}
A_n&:= \frac{\Delta_n}{\tau_{k}(v\Delta_n)^2} \sum_{i=vk}^{[t/\Delta_n]} \left( 
(\Delta_{i,k}^{n,v} X)^2  - \Big( \sum_{j=0}^l \sigma_{(i-vk)\Delta_n-\theta_j} \Delta_{i,k}^{n,v} G^{(j)} \Big)^2 \right) \\[1.5 ex]
B_{n,m} &:= \frac{\Delta_n}{\tau_{k}(v\Delta_n)^2} \left( \sum_{i=vk}^{[t/\Delta_n]} 
\Big( \sum_{j=0}^l \sigma_{(i-vk)\Delta_n-\theta_j} \Delta_{i,k}^{n,v} G^{(j)} \Big)^2 
- \sum_{r=1}^{[mt]} \sum_{i\in I_m(r)}   \Big( \sum_{j=0}^l \sigma_{(r-1)/m-\theta_j} \Delta_{i,k}^{n,v} G^{(j)} \Big)^2 \right)
\\[1.5 ex]
C_{n,m}&:= \frac{\Delta_n}{\tau_{k}(v\Delta_n)^2}
\sum_{r=1}^{[mt]} \sum_{i\in I_m(r)}   \Big( \sum_{j=0}^l \sigma_{(r-1)/m-\theta_j} \Delta_{i,k}^{n,v} G^{(j)} \Big)^2 
- \frac 1m \sum_{r=1}^{[mt]} \int_0^\infty   \sigma_{(r-1)/m -\theta}^2  \pi_k (d\theta)
\\[1.5 ex]
D_m &:= \frac 1m \sum_{r=1}^{[mt]} \int_0^\infty   \sigma_{(r-1)/m -\theta}^2  \pi_k (d\theta)- 
\int_0^\infty \left( \int_{-\theta}^{t-\theta } \sigma_s^2 ds\right) \pi_k (d\theta)
\end{align*}
with
\[
I_m(r):=\Big\{i|~ i\Delta_n\in \Big( \frac{r-1}{m}, \frac rm \Big] \Big\}.
\]
Let us give an interpretation to the introduced decomposition. The term $A_n$ is the error associated with the crucial 
approximation introduced in \eqref{xapprox}. In a second step we divide the interval $[0,t]$ into $[mt]$ sub-blocks
and freeze the intermittency $\sigma$ in the beginning of each block; the associated error is represented by $B_{n,m}$.
Within each sub-block we apply the law of large numbers to the Gaussian part. The error of this procedure is denoted by
$C_{n,m}$. Finally, $D_m$ represents the error of a Riemann sum approximation. Next we will prove that 
\[
\lim_{m\rightarrow \infty} \limsup_{n\rightarrow \infty} \mathbb P(|A_n+B_{n,m} +C_{n,m} +D_m|>\epsilon)=0,
\]
for any $\epsilon >0$. This will complete the proof of Theorem \ref{th1}.    \\ \\
\textit{The term $A_n$.} The convergence $A_n \toop 0$ is shown exactly as in \cite[Section 7.3]{BCP11}. Therein
the proof is given for the case of a single singularity at $0$. However, it directly extends to the case
of multiple singularities. \qed \\ \\
\textit{The term $B_{n,m}$.} Observe that 
\begin{align*}
|B_{n,m}|&\leq \frac{\Delta_n}{\tau_{k}(v\Delta_n)^2} 
\sum_{r=1}^{[mt]}
\sum_{j_1,j_2=0}^l  \sup_{s\in \Big(\frac{r-2}{m}, \frac rm \Big]}
|\sigma_{(r-1)/m-\theta_{j_1}} \sigma_{(r-1)/m-\theta_{j_2}} - \sigma_{s-\theta_{j_1}} \sigma_{s-\theta_{j_2}}| \\[1.5 ex]
&\times \Big| \sum_{i\in I_m(r)}    \Delta_{i,k}^{n,v} G^{(j_1)} \Delta_{i,k}^{n,v} G^{(j_2)} \Big| + R_{n,m}
\end{align*}
with $\lim_{m\rightarrow \infty} \limsup_{n\rightarrow \infty} \mathbb P(|R_{n,m}|>\epsilon)=0$. The dominating term converges
in probability to 
\[
B_m:= \frac 1m \sum_{r=1}^{[mt]} \sum_{j\in \mathcal A} \pi (\theta_j)\sup_{s\in \Big(\frac{r-2}{m}, \frac rm \Big]}
|\sigma_{(r-1)/m-\theta_{j}}^2 - \sigma_{s-\theta_{j}}^2 |
\]
as $n\rightarrow \infty$ due to Theorem \ref{th5}(ii) and convergence
$\tau_{k,j}(v\Delta_n)^2/\tau_{k}(v\Delta_n)^2\rightarrow \pi (\theta_j )$. 
Since $\sigma$ is c\'adl\'ag and bounded on compact intervals, we conclude
that 
\[
B_m \toop 0 \qquad \text{as } m\rightarrow \infty, 
\]
which completes the proof of this part. \qed \\ \\
\textit{The term $C_{n,m}$.} According to Theorem \ref{th5}(ii), we have
\[
\Delta_n \sum_{i\in I_m(r)} \frac{\Delta_{i,k}^{n,v} G^{(j_1)} \Delta_{i,k}^{n,v} G^{(j_2)}}
{\tau_{k,j_1}(v\Delta_n) \tau_{k,j_2}(v\Delta_n)} \toop \delta_{j_1,j_2} m^{-1}t
\]
for any $r$. Furthermore, the proof of Proposition \ref{prop1} shows that 
\[
\frac{\tau_{k,j}(v\Delta_n)^2}{\tau_{k}(v\Delta_n)^2} \rightarrow \pi (\theta_j ).
\]
Thus, for any fixed $m$, we conclude that 
\[
C_{n,m} \toop 0 \qquad \text{as } n\rightarrow \infty,
\]
which completes the proof of this part. \qed \\ \\
\textit{The term $D_{m}$.} Recall that the measure $\pi$ is finite and the process $\sigma$ is c\'adl\'ag bounded.
Hence, the convergence $D_m\toop 0$ as $m\rightarrow \infty$ follows by Lebesgue integrability. 
This completes the proof of Theorem \ref{th1}. \qed \\ \\   
\textit{Proof of Theorem \ref{th3}(i)}: The proof of this result follows along the same lines as the previous one. 
For the treatment of the terms $B_{n,m}$ and $C_{n,m}$ we use the convergence in probability
\begin{align*}
&\Delta_n \sum_{i=vk}^{[t/\Delta_n]} \frac{(\Delta_{i,k}^{n,v} G^{(j_1)} \Delta_{i,k}^{n,v} G^{(j_2)})^2}
{(\tau_{k,j_1}(v\Delta_n) \tau_{k,j_2}(v\Delta_n))^2} \ucp t \qquad \text{when } j_1 \not = j_2, \\[1.5 ex]
& \Delta_n \sum_{i=vk}^{[t/\Delta_n]} \frac{(\Delta_{i,k}^{n,v} G^{(j)} )^4}
{(\tau_{k,j}(v\Delta_n) )^4} \ucp 3t,
\end{align*}
which follows from a general result of \cite[Theorem 1]{BCP11}. The remaining proof of Theorem \ref{th1}
applies directly to  Theorem \ref{th3}(i). \qed

\subsection{Proof of Theorem \ref{th2}}
Here we use a slightly different decomposition than in the proof of Theorem \ref{th1}. Observe that 
\begin{align} \label{cltdec}
\Delta_n^{-1/2} \left(\frac{\Delta_n}{\tau_{k}(v\Delta_n)^2} QV(X,k,v\Delta_n)_t -   QV(X,k)_t \right)
= \widetilde A_n^{v} + \widetilde B_{n,m}^{v} + \widetilde C_{n,m}^{v} + \widetilde D_{n}^{v},
\end{align}
where
\begin{align*}
\widetilde A_n^{v}&:= \frac{\Delta_n^{1/2}}{\tau_{k}(v\Delta_n)^2} \sum_{i=vk}^{[t/\Delta_n]} \left( 
(\Delta_{i,k}^{n,v} X)^2  - \Big( \sum_{j=0}^l \sigma_{(i-vk)\Delta_n-\theta_j} \Delta_{i,k}^{n,v} G^{(j)} \Big)^2 \right) \\[1.5 ex]
\widetilde B_{n,m}^{v}  &:= 
\Delta_n^{1/2} \left( \frac{1}{\tau_{k}(v\Delta_n)^2} \sum_{i=vk}^{[t/\Delta_n]} 
\Big( \sum_{j=0}^l \sigma_{(i-vk)\Delta_n-\theta_j} \Delta_{i,k}^{n,v} G^{(j)} \Big)^2 
- \Delta_n \sum_{i=vk}^{[t/\Delta_n]} \int_0^\infty   \sigma_{(i-vk)\Delta_n- \theta}^2  \pi_k (d\theta) \right) \\[1.5 ex]
&- \Delta_n^{1/2} \left( \frac{1}{\tau_{k}(v\Delta_n)^2}
\sum_{r=1}^{[mt]} \sum_{i\in I_m(r)}   \Big( \sum_{j=0}^l \sigma_{(r-1)/m-\theta_j} \Delta_{i,k}^{n,v} G^{(j)} \Big)^2 
\right. \\[1.5 ex]
&-\left.   \frac 1m \sum_{r=1}^{[mt]} \int_0^\infty   \sigma_{(r-1)/m -\theta}^2  \pi_k (d\theta) \right) \\[1.5 ex]
\widetilde C_{n,m}^{v}&:= \Delta_n^{1/2} \left( \frac{1}{\tau_{k}(v\Delta_n)^2}
\sum_{r=1}^{[mt]} \sum_{i\in I_m(r)}   \Big( \sum_{j=0}^l \sigma_{(r-1)/m-\theta_j} \Delta_{i,k}^{n,v} G^{(j)} \Big)^2 
-\frac 1m \sum_{r=1}^{[mt]} \int_0^\infty   \sigma_{(r-1)/m -\theta}^2  \pi_k (d\theta) \right) \\[1.5 ex]
\widetilde D_{n}^{v} &:= \Delta_n^{1/2} \left( 
\Delta_n \sum_{i=vk}^{[t/\Delta_n]} \int_0^\infty   \sigma_{(i-vk)\Delta_n- \theta}^2  \pi_k (d\theta)
- \int_0^\infty \left( \int_{-\theta}^{t-\theta } \sigma_s^2 ds\right) \pi_k (d\theta) \right)
\end{align*}
In the next step we will prove that 
\[
\lim_{m\rightarrow \infty} \limsup_{n\rightarrow \infty} \mathbb P(\|\widetilde A_n^v+ \widetilde B_{n,m}^v  +
\widetilde D_n^v\|_{\infty}>\epsilon)=0,
\]
for any $\epsilon >0$. The term $\widetilde C_{n,m}^{v}$ will give us the central limit theorem. More precisely,
we will show that
\[
(\widetilde C_{n,m}^{1}, \widetilde C_{n,m}^{2}) \stab C_m \qquad \text{as } n\rightarrow \infty,  
\]
for some process $C_m$ and $C_m \stab L$, where the process $L$ is defined at \eqref{CLT}. This would
complete the proof of Theorem \ref{th2}. \\ \\
\textit{Term $\widetilde A_n^{v}$.} The convergence $\widetilde A_n^{v} \ucp 0$ has been proved for $k=1$
in \cite[Section 7]{BCP11} and for $k=2$ in \cite[Section 5.2]{BCP13} (the latter proof easily extends to any
$k\geq 2$). Although both results are only valid for the case of single singularity at $0$, they extend
to the case of multiple singularities exactly as in the proof of Theorem \ref{th1}. \qed \\ \\
\textit{Term $\widetilde B_{n,m}^v$.} The negligibility of the quantity $\widetilde B_{n,m}^v$ is proven 
by means of fractional calculus in a recent work \cite[Section 4]{CNP14}. \qed \\ \\
\textit{Term $\widetilde C_{n,m}^v$.} We recall first that
\[
\frac{\tau_{k,j}(v\Delta_n)^2}{\tau_{k}(v\Delta_n)^2} = \pi (\theta_j ) + o(\Delta_n^{1/2}) \quad j\in \mathcal A,
\qquad \frac{\tau_{k,j}(v\Delta_n)^2}{\tau_{k}(v\Delta_n)^2} = o(\Delta_n^{1/2}) \quad j \not \in \mathcal A,
\]
which follows from the proof of Proposition \ref{prop1} and condition \eqref{robustness}. 
Define the statistics
\begin{align*}
S_{n,m}^{v,j_1,j_2} (r):=
\Delta_n^{1/2} \left( \frac{1}{\tau_{k}(v\Delta_n)^2} \sum_{i\in I_m(r)} \Delta_{i,k}^{n,v} G^{(j_1)} \Delta_{i,k}^{n,v} G^{(j_2)} 
-\delta_{j_1,j_2} \frac{\pi(\theta_{j_1})}{m} \right), \qquad j_1,j_2\in \mathcal A.
\end{align*}
Applying Theorem \ref{th5}(ii) and the properties of stable convergence, we conclude that
\begin{align*}
&\left(\sigma_{(r-1)/m-\theta_j}, S_{n,m}^{v,j_1,j_2} (r) \right)_{j,j_1,j_2\in \mathcal A}^{r=1,\ldots, m, ~v=1,2} \\[1.5 ex]
&\stab \left(\sigma_{(r-1)/m-\theta_j}, 
\sqrt{\pi(\theta_{j_1})\pi(\theta_{j_2})} \Big(V_{k,v}^{j_1,j_2} \Big(\frac{r}{m} \Big) - V_{k,v}^{j_1,j_2} \Big(\frac{r-1}{m}\Big) \Big) 
\right)_{j,j_1,j_2\in \mathcal A}^{r=1,\ldots, m,~v=1,2},
\end{align*} 
where the process $V$ is defined at \eqref{qvgclt}. Next, we observe that
\begin{align*}
\Delta_n^{1/2} \frac{1}{\tau_{k}(v\Delta_n)^2} \sum_{i\in I_m(r)} \Delta_{i,k}^{n,v} G^{(j_1)} \Delta_{i,k}^{n,v} G^{(j_2)}
=o_{\mathbb P}(1)
\end{align*}
when $j_1\not \in \mathcal A$ or $j_2\not \in \mathcal A$, which is due to \eqref{neglig}.
Hence, it holds that
\begin{align*}
\widetilde C_{n,m}^{v}=  
\sum_{r=1}^{[mt]}   \sum_{j_1,j_2\in \mathcal A} \sigma_{(r-1)/m-\theta_{j_1}} 
\sigma_{(r-1)/m-\theta_{j_2}} S_{n,m}^{v,j_1,j_2} (r)+ o_{\mathbb P}(1).
\end{align*}
Now, applying the continuous mapping theorem for stable converge and recalling the identity \eqref{lambdaid}, 
we deduce that
\[
(\widetilde C_{n,m}^{1}, \widetilde C_{n,m}^{2}) \stab \sum_{r=1}^{[mt]} 
\left( \int_0^{\infty} \sigma_{\frac{r-1}{m}-\theta}^2 ~\pi_k(d\theta) \right) \Lambda^{1/2}_k (B_{\frac{r}{m}} - B_{\frac{r-1}{m}})
\qquad \text{as } n\rightarrow \infty, 
\] 
where $\Lambda_k$ and $B$ are defined at \eqref{CLT}. Finally,
\[
\sum_{r=1}^{[mt]} 
\left( \int_0^{\infty} \sigma_{\frac{r-1}{m}-\theta}^2 ~\pi_k(d\theta) \right) \Lambda^{1/2}_k (B_{\frac{r}{m}} - B_{\frac{r-1}{m}})
\stab L_t \qquad \text{as } m\rightarrow \infty, 
\] 
which completes this step. \qed \\ \\
\textit{Term $\widetilde D_n^v$.} Since $\sigma$ is H\"older continuous of order $\gamma$ with $\gamma >1/2$, 
we readily deduce that $\widetilde D_n^v \ucp 0$. This completes the proof of Theorem
\ref{th2}. \qed

\end{document}